\documentclass[10pt]{svmult}

\usepackage{graphicx,amssymb,amsmath}
\usepackage{color}
\usepackage{setspace}
\usepackage{CJKutf8} 
\usepackage[T1]{fontenc}
\usepackage{ucs} 
\usepackage[encapsulated]{CJK} 
\def\R{\mathbb R}
\def\C{\mathbb C}

\def\N{\mathbb N}
\def\T{\mathbb T}

\let\o\omega
\def\e{\boldsymbol\varepsilon}

\let\a\alpha
\let\cal\mathcal
\let\hat\widehat

\def\mfa{{\mathfrak a}}
\let\l\lambda
\newtheorem{thm}{Theorem}
\theoremstyle{definition}
\newtheorem{defn}[thm]{Definition}

\usepackage[utf8]{inputenc} %latin1
\title*{An overview of dualities in non-commutative harmonic analysis}

\renewcommand\thesection{\arabic{section}}

\author{ Yulia Kuznetsova}
\institute%[Besan\c con] % (optional, but mostly needed)
{ University of Franche-Comt\'e, Besan\c con, France}

\begin{document}
%{\footnotesize
%\def\o{\boldsymbol\omega}

%\title*{Contribution Title}
%\author{Name of First Author and Name of Second Author} %order: first (given) name family name (no abbreviations!)
%\institute{Name of First Author \at Name, Address of Institute \at \email{name@email.address} % %order: first (given) name family name (first name abbreviated!)
%\and Name of Second Author \at Name, Address of Institute }  %Provide only (preferably university) e-mail of corresponding author! 

\maketitle
 
% \tableofcontents

\section{Introduction}

The aim of this mini-course is to give an overview of tools of the non-commutative harmonic analysis. By volume limitations, it cannot be complete, but tries to be as wide as possible and accessible to starting PhD students.

Most of the general theory presented below can be bound in the books of Jacques Dixmier \cite{dixmier} and Gerald Folland \cite{folland}. On spherical theory, the recommended reference is Sigurdur Helgason, Groups and Geometric Analysis \cite{helgason}. And on quantum groups, Thomas Timmermann, An Invitation to Quantum Groups and Duality \cite{timm}. Other sources are research articles.

\emph{\textbf{Duality and Fourier transform}}

Duality is one of the main tools in classical harmonic analysis. The Fourier transform $\cal F$ takes convolution to pointwise multiplication, thus allowing to simplify differential equations and to estimate their solutions. Keep in mind that solutions of differential equations are convolutions with Green functions, so we are particularly interested in convolution operators $T: f \mapsto f*k$; the kernel $k$ can be a function, measure, or a distribution.

This theory is fully available on an Abelian locally compact group $G$: the group of characters of $G$
$$
\hat G = \{ \phi: G\to \T {\rm \ continuous \ homomorphism } \}
$$
is also an Abelian LCG, and $\hat{\hat G} \simeq G$.
The Fourier transform $\cal F$ satisfies $\cal F(f*g) = \cal F(f) \cal F(g)$ and is isometric on $L^2$, thus it is easy to say when $T$ is bounded on $L^2$: exactly when $\cal F(k)$ is bounded.

If $G$ is non-abelian, such a duality is not possible. The first natural analogue of the group of characters is $\hat G$ : the set of all unitary irreducible representations of $G$.
But $\hat G$ is no more a group. Three solutions exist:
\begin{itemize}
\item Work with $\hat G$ as a topological or a measured space
\\ $\longrightarrow$ Plancherel measure, leads to an isometry on $L^2$ in many cases
\item Consider only finite-dimensional representations
\\ $\longrightarrow$ tensor categories, especially useful for compact groups
\item Study group algebras
\\ $\longrightarrow$ Duality of quantum groups
\end{itemize}
And a special domain is the study of subclasses of operators which admit Fourier-like methods (analysis on symmetric spaces).

I will speak of all of them except for the second one which could not enter in the volume available. A good reference for this approach is \cite{neshveyev}.

\section{Unitary dual}

Abelian groups are well described by their characters having scalar values, and the Fourier transform of a function is again a scalar-valued function. In the non-commutative case, representations act on Hilbert spaces, possibly of infinite dimension, and the Fourier transform will be an operator-valued function on $\hat G$.

Recall that every unitary representation $\pi$ of $G$ on a Hilbert space $H_\pi$ generates a representation of $L^1(G)$ by
\begin{equation}\label{pi(f)}
\pi(f)  = \int_G f(g) \pi(g) dg;
\end{equation}
the integral is understood as an operator such that for every $x,y\in H_\pi$
$$
\langle \pi(f)x,y\rangle = \int_G f(g) \langle \pi(g)x,y\rangle dg.
$$
This representation extends by continuity to $C^*(G)$.
Conversely, every *-representation of $C^*(G)$ defines a unitary representation of $G$
(provided that $\pi(C^*(G))H$ is dense in $H$).
A representation $\pi$ is called irreducible if it has no nontrivial closed invariant subspaces.
If $G$ is abelian, every irreducible $\pi$ is one-dimensional.

Representations $\pi$, $\rho$ of $G$ on Hilbert spaces $H_\pi$, $H_\rho$ are called equivalent if there exists a unitary $T:H_\pi\to H_\rho$ such that  $T\pi(g) = \rho(g)T$ for every $g\in G$.

\begin{defn}
$\hat G$ is the set of equivalence classes of {\it unitary irreducible} representations.
\end{defn}

\def\mfg{\mathfrak{g}}

\smallskip
\emph{\textbf{Nilpotent case}}

It is usually a complicated task to describe $\hat G$. For a nilpotent simply connected Lie group $G$, a powerful orbit method is available.
Examples of such groups include strictly upper-triangular matrices and Heisenberg groups $\R\ltimes\R^n$.

Let $\mfg$ denote the Lie algebra of $G$. At $X\in \mfg$, the {\it adjoint action} of $g\in G$ is
$$
Ad (g): \exp(tX) \mapsto g\exp(tX)g^{-1}.
$$
On $\mfg^*$, this generates the {\it coadjoint action}:
$$
Ad^*(g): \l \mapsto \l \circ Ad (g).
$$
We say that $\l\sim \mu$ if $\l = Ad^*(g)\mu$ for some $g$. This is an equivalence relation. Without going into further detail, let us state
\begin{theorem}[Kirillov \cite{kirillov}]
$\hat G$ is in bijection with $\mfg^* / \sim$.
\end{theorem}

\smallskip
\emph{\textbf{Heisenberg group}}

\def\H{\mathbb{H}}
This example will help to see the identification above. We set $\H = \R \ltimes \R^2$, or the group of matrices
$
\begin{pmatrix} 1 & a& c\\
0 & 1& b\\
0& 0& 1
\end{pmatrix}
$.
\iffalse{
with
$$
\begin{pmatrix} 1 & a_1& c_1\\
0 & 1& b_1\\
0& 0& 1
\end{pmatrix}
\begin{pmatrix} 1 & a_2& c_2\\
0 & 1& b_2\\
0& 0& 1
\end{pmatrix}
=
\begin{pmatrix} 1 & a_1 + a_2 & c_1 + c_2 + a_1 b_2\\
0 & 1& b_1 + b_2\\
0& 0& 1
\end{pmatrix}
$$
}\fi

Its Lie algebra is $\mfg = \langle X,Y,Z\rangle$ with
$$
X = \begin{pmatrix} 0 & 1& 0\\
0 & 0& 0\\
0& 0& 0
\end{pmatrix} \
Y = \begin{pmatrix} 0 & 0& 0\\
0 & 0& 1\\
0& 0& 0
\end{pmatrix} \
Z = \begin{pmatrix} 0 & 0& 1\\
0 & 0& 0\\
0& 0& 0
\end{pmatrix},
$$
so that $[X,Y]=Z$ and $[X,Z]=[Y,Z] = 0$.

The space $\mfg^*$ is 3-dimensional and the coadjoint action (we do not calculate it here) is 
$$
Ad^*(a,b,c) : (\mu,\nu,\l) \mapsto (\mu+b\l, \nu-a\l, \l).
$$
We see that the orbit of $(\mu,\nu,0)$ is only this point, and if $\l\ne0$, the orbit is the plane 
$\Omega_\l = \{ (\mu,\nu,\l): (\mu,\nu)\in \R^2\}$. Denote $\Omega = \{ \Omega_\l: \l\in\R\}$.

For $\mu,\nu\in\R$, the one-dimensional representation corresponding to $(\mu,\nu,0)$ is given by
$$
\pi_{\mu,\nu}(a,b,c) = e^{ 2\pi i (a\mu + b\nu) };
$$
for $\l\in\R$, infinite-dimensional irreducibles act on $f\in L^2(\R)$ as
$$
\pi_\l (a,b,c) (f) (t) = e^{ 2\pi i \l (c+bt) } f(t+a).
$$
%Kirillov p.67

By direct calculations, one can verify that for $f\in L^1(\H) \cap L^2(\H)$,
\iffalse{
$$
\pi_{x,y}(f) = \int_H f(a,b,c) e^{ 2\pi i(ax+by)} dadbdc = \hat f(x,y,0);
$$
$$
\pi_z(f)
$$
}\fi
$$
\| f\|_2 ^2 = \int_\R Tr [\pi_\l(f)^* \pi_\l(f) ] \,|\l| d\l.
$$
The measure
\begin{equation}\label{plancherel-heisenberg}
\mu(\{\pi_{\mu,\nu}\})=0, \quad \mu|_{\Omega} = |\l| d\l
\end{equation}
is called the Plancherel measure on $\hat {\H}$, and we will discuss it below in the general case.
% Folland p.239
% Theorem: p.234
% Dixmier [29, 18.8]

\smallskip
\emph{\textbf{Topology on the unitary dual}}

This dual is more often considered as a measure space, but one starts with a topology. For $\pi\in \hat G$, its kernel $\ker\pi$ is an ideal in $C^*(G)$.
We define the {\it closure} of any $E\subset \hat G$:
$$
\overline{E} = \{ \pi\in \hat G: \pi \equiv0 {\rm\ on\ } I_E = \cap_{\rho\in E} \ker\rho \}.
$$
There exists (by Kuratowski theorem) a topology on $\hat G$ in which this is the closure map.
One calls it the Fell, or hull-kernel topology. $\hat G$ is not always Hausdorff.

\smallskip
If $G$ is {\it abelian}, then
$C^*(G) \simeq C_0(\hat G)$.
For $\rho\in \hat G$, we have $\ker \rho = \{ f\in C_0(\hat G): f(\rho)=0 \}$,
so that $I_E = \{ f\in C_0(\hat G): f|_E \equiv 0 \}$. Then
$$
\overline{E} = \{ \pi\in \hat G: f(\pi)=0 {\rm \ for\ every\ }f\in C_0(\hat G) {\rm\ vanishing\ on\ } E \}.
$$
Clearly this is the usual closure of $E$.

\subsection{Plancherel measure}

For $f\in L^1(G)$ and $\pi\in\hat G$, set $\hat f(\pi) = \pi(f)$, as given by the formula \eqref{pi(f)}: this defines the noncommutative Fourier transform.
Sometimes one sets instead $\pi(f)  = \int_G f(g) \pi(g^{-1}) dg$.

For groups of so called ``Type I'', there is a direct analogy with the classical case:
\begin{theorem}
Let $G$ be a \underline{unimodular Type I} group.
There exists a unique positive measure $\mu$ on $\hat G$ such that for every $f\in L^1(G)\cap L^2(G)$
\begin{equation}\label{Fourier-L2}
\| f\|_2 ^2 = \int_{\hat G} Tr \big( \hat f(\pi)^* \hat f(\pi) \big) d\mu(\pi).
\end{equation}
\end{theorem}
There exists a non-unimodular version \cite[Theorem 7.50]{folland}
but the group should be ``almost'' Type I.

\medskip
 %Dixmier 9.1
Let us now define the term ``Type I''. Recall that a subalgebra $A$ in $B(H)$ (the space of bounded operators on a Hilbert space $H$) is called a von Neumann algebra if it is closed in strong or, equivalently, weak operator topology. This happens if and only if $A$ is equal to its second commutant $A''$; and for a set $S\subset B(H)$, its commutant is defined as $S' = \{ y\in B(H): yx=xy\ \forall x\in S\}$. Next, a von Neumann algebra $A$ is a factor if $A\cap A'=\C 1$.

Now, being Type I means: if $\pi(G)''$ is a factor, then it is of type I, that is, isomorphic to $B(H)$ (of finite or infinite dimension).This is also equivalent to the following: if $\pi$ is irreducible, then $\pi\big( C^*(G) \big)$ contains all compact operators. If $G$ is not of Type I, then it has necessarily both Type II and Type III representations \cite{glimm}. An example of a non-Type I group is the free group of two generators.

Other terms for this class are GCR and postliminal groups. These notions have appeared independently and have different definitions, but have turned to be the same classes of groups as Type I.
% not necessarily second countable

In relation to the dual, the following facts are known. For any $G$, if $\pi_1 \!=\! \pi_2\in \hat G$, then clearly $\ker \pi_1 \!=\!  \ker \pi_2$.
For Type I groups, the converse is true. Moreover, a second countable group $G$ is Type I if and only if the map $\pi \mapsto \ker\pi$ is injective, and if and only if the topology of $\hat G$ is $T_0$.

Finally, if $G$ is Type I and $\pi$ irreducible, then the von Neumann algebra $\pi(G)''$ admits a trace.

\begin{itemize}
\item Every connected semisimple and nilpotent Lie group is Type I.
\item A discrete group is Type I if and only if it is almost abelian (that is, has an abelian subgroup of finite index).
\end{itemize}

\smallskip
\emph{\textbf{Compact groups}}

In this case, $\hat G$ is discrete: if $\chi_\pi: g\mapsto Tr \,\pi(g)$ is the character of $\pi$, then
$\rho( \chi_\pi) = 0$ if $\pi\ne\rho$ and $\pi(\chi_\pi)\ne0$,
so that $\hat G \setminus\{\pi\}$ is closed, thus $\{\pi\}$ open.

The Plancherel theorem takes the following form:
$$
 \| f\|_2 ^2 = \sum_{\pi\in \hat G} d_\pi Tr \big( \pi(f) \pi(f)^* \big).
$$

\smallskip
\emph{\textbf{Inverse Fourier transform}}

The Plancherel identity implies that for $f,g\in L^1\cap L^2(G)$
\begin{align*}
\langle f,g \rangle &= \int_{\hat G} Tr \big( \pi(g)^* \pi(f) \big) d\mu(\pi)
\\&= \int_{\hat G} \int_G Tr \big( \overline{g(x)} \pi(x)^* \pi(f) \big) dx d\mu(\pi),
\end{align*}
so that
$$
f(x) = \int_{\hat G} Tr \big( \pi(x)^* \pi(f) \big) d\mu(\pi).
$$
This allows in particular to define pseudo-differential operators on such groups \cite{ruzh-pseudo}.

\smallskip
\emph{\textbf{Plancherel measure: the support}}

The support of $\mu$ (intersection of all closed sets of full measure)
 is the set of $\rho\in \hat G$ {\it weakly contained} in $\lambda$, that is, such that
%Dixmier 18.8.4 & 18.1.7
for every $\xi\in H_\rho$, $\e>0$, and compact $K\subset G$
 there exist $f_k\in L^2(G)$ with
$$
\sup_{g\in K} \big| \langle \rho(g) \xi,\xi\rangle - \sum_{k=1}^n \tilde f_k*f_k (g) | <\e.
$$
Here $\tilde f(x) = \bar f(x^{-1})$.

$G$ is amenable iff ${\rm supp}\,\mu = \hat G$.
%Folland p.235
This is not the case, for example, for $SL_2(\R)$.
% or the 'ax+b' group.

\smallskip
\emph{\textbf{Example: Heisenberg group}}

Recall that $\hat {\H} = \{ \pi_{\mu,\nu} \} \cup \{ \pi_\l \}$, and only $\pi_\l$ enter in the Plancherel formula \eqref{plancherel-heisenberg}.
For $k\in L^1(\H) \cap L^2(\H)$,
\begin{align*}
\pi_\l (k) f (t) &= \int_H k(a,b,c) e^{ 2\pi i \l (c+bt) } f(t+a) dadbdc
\\ &= \int_\R \cal F_2 \cal F_3 (k) (x-t, -\l, - \l t ) f(x) dx
\end{align*}
so that $\hat k(\pi_\l)$ is an integral operator with the kernel as above.

\smallskip
\emph{\textbf{Plancherel measure: non-unimodular case}}

The Plancherel measure still exists, with Formula \eqref{Fourier-L2} and \eqref{disint-regular} valid, but the group should be ``almost Type I'' \cite[Theorem 7.50]{folland} and we have to adjust the definition of the Fourier transform setting
$$
\hat f(\pi) = \pi(f) D_\pi,
$$
with $\pi(f)$ as in \eqref{pi(f)}.
Here $D_\pi$ are self-adjoint operators on $H_\pi$ such that $D_\pi \pi(g) = \Delta(g)^{1/2} \pi(g) D_\pi$ for every $g\in G$. Their existence is not obvious, but becomes clearer knowing that every $\pi$ can be realized on the space $H_\pi = L^2(G, V_\pi)$ where $V_\pi$ is also a Hilbert space; one sets then $D_\pi f(g) = \Delta(g)^{-1/2} f(g)$ for $f\in H_\pi$ and $g\in G$.

The Plancherel identity for $f,g\in L^1\cap L^2(G)$ becomes
\begin{align*}
\langle f,g \rangle &= %\int_{\hat G} Tr \big( \hat g(\pi)^* \hat f(\pi) \big) d\mu(\pi) \\&= 
\int_{\hat G} \int_G Tr \big( D_\pi \overline{g(x)} \pi(x)^* \pi (f) D_\pi \big) dx d\mu(\pi),
\end{align*}
and the inverse Fourier transform
$$
f(x) = \int_{\hat G} Tr \big( \pi(x)^* \pi(f) D_\pi^2 \big) d\mu(\pi).
$$

\smallskip
\emph{\textbf{`ax+b' group}}

$G$ is the group of matrices
\vskip-8pt
$$
\begin{pmatrix}
a& b\\
0 & 1
\end{pmatrix}
$$
with $a>0$, $b\in\R$, which correspond to the affine maps on the real line: $ g_{a,b} x = ax+b$.

We could consider also $G_n = \{ (a,b): a\in \R, b\in \R^n\}$, acting the same way.
$G$ is non-unimodular: the left Haar measure is $\dfrac1{a^2} da db$, the right one $\dfrac1a da db$,
the modular function is $a^{-1}$.

One can describe $\hat G$ by the orbit method, arriving at: $\hat G = \{ \pi_\l: \l\in\R\} \cup \{\pi^+, \pi^-\}$, acting as
$$
\pi_\l (a,b) = a^{i\l};
$$
$\pi^+$, $\pi^-$ act on $H^+ = L^2(0,+\infty)$ and $H^- = L^2(-\infty,0)$ respectively by the same formula:
$$
\pi^\pm(a,b) \phi(t) = \sqrt a \, e^{ 2\pi ib t} \phi(at). %, \quad \pi^-(a,b) f(t) = \sqrt a \, e^{ 2\pi ib t} f(at).
$$
The Plancherel measure $\mu$ on $\hat G$ is such that
$\mu(\{\pi_\l\})=0$ for every $\l$, and
$$
\| f\|_2^2 = Tr \big[ \hat f(\pi^+)^* \hat f(\pi^+) \big] + Tr \big[ \hat f(\pi^-)^* \hat f(\pi^-) \big]
$$
for $f\in L^2(G)$. The adjustment operators are
$
D^\pm \phi(s) = \sqrt{|s|} \phi(s)
$, and $ \hat f(\pi^\pm) = \pi^\pm(f) D^\pm.
$
One verifies directly that $\hat f(\pi^+)$ is an integral operator with the kernel
\begin{align*}
k_f^+(s,t) % &= \frac{ \sqrt s} t \int_\R f(s^{-1} t ,b) \, e^{2\pi i bs} db\\
&=  \frac{ \sqrt s} t \,\cal F_2 f(s^{-1} t , -s ),
\end{align*}
the Fourier transform being applied in the second variable only. A similar formula holds for $\hat f(\pi^-)$.

\subsection{Disintegration}

The Plancherel measure has much more properties than just the isometry of the Fourier transform on $L^2(G)$. The next property to discuss is the decomposition of representations into direct integrals of irreducibles \cite{dixmier}.

Let $(\Omega,\mu)$ be a measured space and $H_\o$ a Hilbert space for every $\o\in\Omega$.
A {\it measurable field} of Hilbert spaces %$(H_\pi)_{\pi\in\hat G}$ 
is a linear subspace
$\Gamma \subset \prod_{\o\in\Omega} H_\o$ such that
\begin{itemize}
\item $\o \mapsto \| x(\o) \|$ is measurable for every $x\in\Gamma$;
\item there exist $(x_n)\subset\Gamma$ such that $\{ x_n(\o)\}_{n\in\N}$ is total in $H_\o$ for every $\o$;
\item if $\o\mapsto \langle x(\o), y(\o)\rangle$ is measurable for every $x\in \Gamma$, then $x\in\Gamma$.
\end{itemize}
We denote
$$
\Gamma := \int^\oplus H_\o d\mu(\o).
$$
If $A(\o) \in B(H_\o)$ for every $\o$ and $(A(\o) x(\o)) \in\Gamma$ for every $x\in\Gamma$,
 then $A = \int^\oplus A(\o) d\mu(\o)$ is said to be a mesurable field of operators on $\Gamma$. It is bounded if and only if
 $$
 \| A\| = {\rm esssup}_{\o\in\Omega} \| A(\o)\| <\infty.
 $$

\smallskip
\emph{\textbf{Disintegration of representations}}

If $\pi_\o$ is a unitary representation of $G$ on $H_\o$ for every $\o\in\Omega$ and
$(\pi_\o(g))$ is a measurable field of operators for every $g\in G$, then
$ \pi = \displaystyle\int^\oplus \pi_\o d\mu(\o) $, the direct integral of $\pi_\o$,
is a representation of $G$ on $\int^\oplus H_\o$. Conversely,
if $G$ is a second countable locally compact group, then every unitary representation $\pi$ of $G$ on a separable $H$ is unitarily equivalent to a direct integral $\displaystyle\int^\oplus \pi_\o d\mu(\o) $ where
$\pi_\o$ are irreducible for almost all $\o$
\cite[8.5.2]{dixmier}, \cite[7.37]{folland}.

\def\F{\mathbb F}
\let\t\theta

This disintegration is not unique, and this can be seen on the example of the most important representation: of the regular one. The left regular representation $\l$ of $G$ acts at $f\in L^2(G)$ as
$$
[\l(g) f] (h) = f( g^{-1} h).
$$

\emph{\textbf{A non-uniqueness example}} \cite{yoshizawa}. Let $\F_2$ be the free group of two generators $a$ and $b$. Then $\l$ has two decompositions into irreducibles
$$
\l = \int_\T U_\t d\t = \int_\T V_\t d\t
$$
such that none of $U_\t, V_{\t'}$ are equivalent.

We introduce briefly these representations. Let $G_a, G_b$ denote all words in $\F_2$ which end by $a,a^{-1}$ or $b,b^{-1}$ respectively, and set $G_a'=G_a\cup\{e\}$, $G_b'=G_b\cup\{e\}$.
For $\t\in\T$, the space $H_\t$ is the closed linear span of the indicator functions $\{\delta_g: g\in G_b'\}$ which are its orthonormal basis; $U_\t$ acts as
\begin{align*}
U_\t (a) \delta_e &= \t \delta_e, \quad U_\t(a) \delta_g = \delta_{ag}, \ g\in G_b; \\
U_\t(b) \delta_g &= \delta_{bg}.
\end{align*}

For $V_\t$, the construction is the same with $a$ and $b$ interchanged.\\
One can show that $U_\t$, $V_{\t'}$ are irreducible and pairwise inequivalent.

\smallskip
\emph{\textbf{Disintergation over the dual}}

To achieve uniqueness, or otherwise stated, to integrate over the dual, one has to return to Type I groups. Multiplicities occur in this process:
\def\n{\mathfrak n}

\begin{theorem}\cite[7.40]{folland}
Let $G$ be a second countable \underline{Type I} group, $\sigma$ a unitary representation of $G$ on a separable $H$.
Then there exist finite pairwise orthogonal measures $\mu_\n$, $\n\in\N\cup\{\infty\}$, on $\hat G$,
such that
$$
\sigma \simeq \int^\oplus \pi d\mu_1(\pi) \oplus 2 \int^\oplus \pi d\mu_2(\pi)
\oplus \dots \oplus |\N| \int^\oplus \pi d\mu_\infty(\pi).
$$
These measures are unique up to equivalence.
\end{theorem}

In particular, for the regular representation 
\begin{equation}\label{disint-regular}
\lambda = \int_{\pi\in \hat G} \pi\otimes I_{\bar H_{\pi}} \, d\mu(\pi)
\end{equation}
with the Plancherel measure $\mu$; the group should be second countable and Type I \cite[7.44, 7.50]{folland}.
%Dixmier 18.7.7
%non-unimodular: Folland 7.50
If $G$ is compact, we recover the familiar decomposition
$$
\lambda = \mathop{\oplus}\limits_{\pi\in \hat G} d_\pi \pi.
$$

\emph{\textbf{Translation-invariant operators}}

Not only $\lambda(g)$ but every translation invariant operator admits such a disintegration. 
Recall that convolution operators are translation invariant; often, converse is true \cite{wendel}, \cite[Corollary 3.2.1]{f-ruzh}.
If $T$ is bounded on $L^2$, then:\\
\hbox{}\hskip80pt $T$ is left invariant $\Leftrightarrow$ $T\in \l(G)''$,\\
where $S''$ is the bicommutant in $B(L^2(G)))$, or in other terms,
the von Neumann algebra generated by $S$.

In the same assumptions on $G$ as above,
for every $T\in \l(G)''$
\begin{equation}\label{inv-op-disintegr}
T = \int^\oplus_{\hat G} T_\pi \otimes I_{\bar H_{\pi}} \,d\mu(\pi).
\end{equation}
Only the support of $\l$ enters in this decomposition.
If $Tf = k*f$, then
$$
T_\pi = \hat k(\pi) = \int_G k(x) \pi(x) dx.
$$

\smallskip
\emph{\textbf{Heisenberg group: sub-Laplacian}}
\smallskip

One can verify that, on $L^2(\R)$, for every $\l\in\R$
$$
\pi_\l(X) = \sqrt{|\l|} \frac{ \partial } {\partial t} , \quad \pi_\l(Y) = i \,\text{sgn}\l \sqrt{|\l|} \,t , \quad \pi_\l(Z) = 2\pi i\l.
$$
%Kirillov p.60, Fischer p.433

If $\cal L = X^2 + Y^2$ is the sub-Laplacian, then, even if it is unbounded, the disintegration \eqref{inv-op-disintegr} holds for it \cite{f-ruzh}, with
\begin{equation}\label{sublapl-fourier}
\cal L_{\pi_\l} f (t) = |\l| \big( f''(t) - t^2 f(t) \big).
\end{equation}

\section{Symmetric spaces}

Let $G$ be a semisimple Lie group, $K$ its maximal compact subgroup. By Iwasawa decomposition, $G=KAN$ with $A\simeq \R^l$ abelian, $N$ nilpotent. Let $\mfa \simeq \R^l$ and $\mfg$ denote the Lie algebras of $A$ and $G$.

The principal series representations $\pi_{\sigma,\l}$ of $G$ are parametrized by irreducibles $\sigma\in \hat K$ and $\l\in \mfa^*$; if $\sigma=1$, they have $K$-invariant vectors and are called spherical. One can realize them on subspaces of $L^2(K)$.

For every $\l\in \mfa^*$, the function $\phi_\l$ on $G$ is the following coefficient of $\pi_{1,\l}$:
$$
\phi_\l(x) = \langle \pi_\l(x) 1, 1\rangle_{L^2(K)}.
$$

For every $x\in G$ one can write uniquely $x = k_x a_x n_x$; set $H(x) = \log a_x$, then
$$
\phi_\l(x) = \int_K e^{ - (i\l+\rho) \big( H(x^{-1} k) \big) } dk,
$$
where $\rho\in \mfa^*$ is defined by the algebraic structure of $\mfg$.
%GV 1.5.10

\smallskip
\emph{\textbf{Spherical transform}}

For $f$ on $G$, we define the spherical transform
$$
\cal H(f) (\l) = \int_G f(x) \phi_\l(x^{-1}) dx,
$$
so that $\cal H(f)$ is a function on $\mfa^*$. It has many similarities to the classical Fourier transform:
For $f,g$ in the space $C^\infty(K\setminus G/K)$ of $K$-biinvariant smooth functions on $G$, we have
\begin{itemize}
\item $\cal H(f*g) = \cal H(f) \cal H(g)$;
\item $\| \cal H(f)\|_2 = \|f\|_2$;
\item inversion formula: $f (x) = \int_{\mfa^*} \cal H(f)(\l) \phi_\l(x) |{\boldsymbol c}(\l)|^{-2} d\l$,
\end{itemize}
where ${\boldsymbol c}$ is the important Harish-Chandra $\boldsymbol c$-function on $G$.

\smallskip
\emph{\textbf{Spherical functions and differential operators}}

Every $\phi_\l$ is an eigenfunction of all differential operators on $C^\infty(G)$ invariant with respect to $G$ on the left and $K$ on the right. In particular, the structure of a Riemannian manifold on $G/K$ gives rise to its Laplace-Beltrami operator $\cal L$ which can be considered as an operator on $C^\infty(G)$ in the above class.
It is convenient to denote $\cal L_\rho = \cal L+|\rho|^2$, then
$$
\cal L_\rho \phi_\l = |\l|^2 \phi_\l.
$$
For every function of $\cal L_\rho$, we have also $F(\cal L_\rho)\phi_\l = F(|\l|^2) \,\phi_\l$. Being left invariant, it is a convolution operator $f\mapsto f*k_F$, and the inversion formula for the spherical transform implies that
$$
k_F(x) = \int_{\mfa^*} F(|\l|^2) \,\phi_\l(x) |{\boldsymbol c}(\l)|^{-2} d\l.
$$
Fourier analysis methods then apply, since this is just an integral over $\mfa^*\simeq \R^l$. The complicated part is the spherical function, for which one seeks simpler approximations; it is known for example, that for ``good'' $x=e^H \in A$ and $\|H\|\to\infty$,
$$
\phi_\l(e^H) \sim \sum_{k=1}^n {\boldsymbol c}(s_k \l) e^{ ( i s_k \l - \rho)(H) },
$$
%GV 4.4.13
where $\{s_k\}$ form a certain group of reflections on $\mfa$, called the Weyl group.

\smallskip
\emph{\textbf{Examples}}

If $G=SL_n(\R)$, then $K=SO_n(\R)$,\\
 $A = \{ \text{ diagonal matrices in $G$ with positive diagonal } \}$,\\
$N = \{ \text{ upper-triangular matrices with diagonal 1 } \}$, and\\
$G/K \simeq AN = \{$ upper-triangular matrices in $G$ with positive diagonal $\}$.

\bigskip

The `ax+b' group 
 can be identified with $G/K$ where $G = SO(2,1)_o$ is the connected component of the identity of $SO(2,1)$, and $K=SO(2)$.

\section{Traces and weights %Non-commutative integration
}\label{sec-traces}

A {\it weight} on a $C^*$ or von Neumann algebra $A$ is a map $\tau: A^+ \to [0,+\infty]$\\
 such that
\begin{itemize}
\item $\tau$ is additive
\item $\tau(tx) = t\,\tau(x)$ for $t\in\R_+$, $x\in A^+$, with $0\cdot(+\infty) = 0$.
\\ \hskip-20pt If in addition $\tau(xx^*) = \tau(x^*x)$ for all $x\in A$, then $\tau$ is a {\it trace}.
\end{itemize}
A trace $\tau$ is {\it faithful} if $\tau(x)=0$ only for $x=0$, and {\it semifinite} if for any $x\in A^+$ there exists $y\in A^+$ such that $ y\le x$ and $\tau(y)<\infty$.

\medskip
$\star$ On $A=B(H)$: the usual trace $Tr$

$\star$ On $A=C_0(G)$: \ $\tau(f) = \int_G f$, so a weight can be viewed as a non-commutative integral

\smallskip
\emph{\textbf{Plancherel weight / trace}}:
On $A=C^*_r(G)$ or $VN(G)$, we can set
$
\tau( L_k ) = k(e)
$
if $k = f^* *f$, or equivalently
$$
\tau( L_f^*L_f ) = \tau( L_{f^* *f} ) = \| f\|_2^2.
$$
If $G$ is unimodular, this is a trace. In general, just a weight.

We can see that the equality \eqref{Fourier-L2} is in fact a disintegration of $\tau(L_f^*L_f)$ into $\tau_\pi\big( \hat f(\pi)^* \hat f(\pi) \big)$ where every $\tau_\pi$ is the usual trace.

\smallskip
\emph{\textbf{Example:}} Consider again the Heisenberg group $\mathbb H$ and the sub-Laplacian $\cal L = X^2 + Y^2$ on it. For any $u\in\R$, by \eqref{Fourier-L2}
$$
\tau\big( I_{(0,u)}(\cal L) \big)
= \int_\R \tau_\l\big( I_{(0,u)}(\cal L_{\pi_\l} ) \big) |\l| d\l;
$$
the operator $\cal L_{\pi_\l}$ is given by \eqref{sublapl-fourier} and its eigenvalues are known, so one can show \cite{ak-ruzh} that
\begin{equation}\label{tau_L}
\tau\big( I_{(0,u)}(\cal L ) \big) = \int_\R \#\{ k: s_{k,\l} < u\} |\l| d\l \lesssim u^{3/2}.
\end{equation}

\smallskip
\emph{\textbf{Non-commutative $L_p$ spaces}}

Let $M$ be a von Neumann algebra with a faithful semifinite normal trace $\tau$ (normal means that $\tau(\sup x_i) = \sup \tau(x_i)$ for every bounded increasing net in $A^+$).
For $x\in M$, define
$$
\| x\|_p = \tau( |x|^p )^{1/p} %= \tau\big( (x^* x) ^{p/2} \big) ^{1/p}.
$$
and $L_p(M)$ as the completion of $\{ x\in M: \|x\|_p<\infty \}$.

For $p=\infty$, $\|\cdot\|_\infty$ is the usual norm of $M = L^\infty(M)$.

\bigskip
$\star$ If $M=L^\infty(X)$ with $\tau(f) = \int_X f d\tau $, then $L^p(M) = L^p(X)$.

$\star$ If $M=\l(G)''=VN(G)$, then with the Plancherel weight, for $f\in L^2(G)$
$$
\| L_f \|_{L^2(M)} = \tau( L_f^* L_f)^{1/2} = \|f\|_2.
$$

\smallskip
\emph{\textbf{Weak $L_p$ spaces}}

One can define {\sl generalized singular numbers} $\mu_t(x)$ for $x\in M$, $t\in\R_+$:
$$
\mu_t(x) = \inf\{ \l>0: \tau\big( I_{( |\l|,+\infty)} (x) \big) <t \}
$$
with the help of spectral projections of $x$, and set
$$
\|x\|_{p,\infty} = \sup_{t>0} t^{1/p} \mu_t(x).
$$
The obtained completed space $L^{p,\infty}(M)$ is called the non-commutative weak $L_p$ space, or Lorentz space.

For a decreasing function $f\ge0$ on $\R_+$, we get $\mu_t(f) = f(t)$. If $x$ is a compact operator on a Hilbert space, then $\mu_n(x) = \l_n(|x|)$ are the eigenvalues of $|x|=(x^*x)^{1/2}$.

\smallskip
\emph{\textbf{Non-unimodular case}}

If $G$ is not unimodular, the Plancherel weight is not a trace.
Non-commutative $L^p$ spaces can be defined still, with complex interpolation available for them; %a construction by Haagerup;
but the Lorentz spaces not, it is an open problem.

There are three definitions of non-tracial $L^p$ spaces, each of them with its advantages and disadvantages: by Haagerup (exposed by Terp \cite{terp}), Izumi \cite{izumi} and Hilsum \cite{hilsum}. We will not go further here.

\smallskip
\bigskip
\emph{\textbf{Quantum groups}}

This notion generalizes both $C_0(G)$ and $C^*_r(G)$, and provides a duality between them. An equivalent theory is built upon von Neumann algebras rather than $C^*$-algebras, and includes $L^\infty(G)$ and $VN(G)$.
Many new examples, not coming from groups, are included in this theory.

\begin{definition} [Kustermans, Vaes \cite{kust-vaes}] A {\it locally compact quantum group} (LCQG) is a von Neumann algebra $M$ with a map $\Delta: M\to M\otimes M$, called comultiplication, and two weights $\phi,\psi$ on $M$ such that
\begin{itemize}
\item $\Delta$ is a unital normal $*$-homomorphism;
\item $\Delta$ is coassociative: $(1\otimes\Delta)\Delta = (\Delta\otimes1)\Delta$;
\item $\phi,\psi$ are normal and semifinite;
\item $\phi$ is left invariant and $\psi$ is right invariant: for $\omega \in M_*^+$ and $x\in M^+$,
\\ if $\phi(x)$ or $\psi(x)$ respectively is finite, then
\begin{align*}
\phi \big( \omega\otimes1)\Delta(x) \big) = \phi(x) \omega(1),\\
\psi \big( 1\otimes\omega)\Delta(x) \big) = \psi(x) \omega(1).
\end{align*}
\end{itemize}
\end{definition}

On $M=L^\infty(G)$, the weights are integrals with respect to the left and right Haar measures. On $M=L^\infty(G)$, they are Plancherel weights.

\smallskip
\emph{\textbf{Duality:}}
 There exists a dual LCQG $(\hat M, \hat\Delta, \hat\phi, \hat\psi)$ such that $\,\widehat{\hat M} = M$. The algebras $L^\infty(G)$ and $VN(G)$ are duals of each other.
There are also analogues of the Fourier transform: $*$-homomorphisms with dense ranges
$$
\l: M_* \to \hat M, \quad \hat\l: \hat M_* \to M.
$$
For $M=L^\infty(G)$, we have $M_* = L^1(G)$ and $\l: f\mapsto L_f \in VN(G)$. The map $\l$ is thus the left regular representation of $L^1(G)$. If~$G$ is abelian, then $VN(G)$ is identified with $L^\infty(\hat G)$ via the usual Fourier transform, with $\hat{L_f g} = \hat{f*g} = \hat f \hat g$, so modulo this identification $L_f$ is the operator of multiplication by $\hat f$ and thus deserves its name as the Fourier transform of $f$.

For $1\le p\le 2$, there are $L_p$-Fourier transforms (Cooney \cite{cooney}, Caspers \cite{caspers}):  $\cal F_p$ maps $L^p(M)$ to $L^{p'}(\hat M)$, with $1/p+1/p'=1$.
In the case of $M=L^\infty(G)$, we have $L^p(M) = L^p(G)$, and $\cal F_p: f \mapsto  L_f \Delta^{1/p'} \in L^{p'}(VN(G))$, where we need to multiply first by a power of the modular function.

$M$ is said of Kac type is $\phi,\psi$ are traces. The example above shows that in the non-Kac case, in particular for a non-unimodular group, $\cal F_p$ are different operators for different $p$, as opposed to the classical case where only their domains are different.

\smallskip
\emph{\textbf{Hausdorff-Young inequality}}

If $G$ is an abelian locally compact group, then for $1\le p\le 2$ % and $1/p+1/p'=1$
$
\| \hat f\|_{p'} \le C_p \|f\|_p
$.
On quantum groups this is also true \cite{caspers}:
$$
\| \cal F_p(f) \|_{L^{p'}(\hat M)} \le \|f\|_{L^p(M)},
$$
with an equality for $p=2$.
In the Kac case, a stronger inequality holds \cite{haonan.zhang}, involving the weak norm $\|f\|_{L^{p,p'}(M)}$.

\section{Multipliers}

In the classical case, a multiplier with a symbol $m$ is the operator $T_m = \cal F^{-1} ( m \hat f)$.
On a quantum group of Kac type, one can adopt the same approach: for $m\in \hat M$, we set
$$
T_m ( y) = \cal F^{-1} \big( m \cal F (y) \big),
$$
$y\in M$. We assume from now that both $M$ and $\hat M$ are of Kac type. Up to my knowledge, there is no theory of non-Kac multipliers, as those would require a coordination of $\cal F_p$ transforms for different $p$. An exception is given by spherical multipliers on symmetric spaces: these groups are non-unimodular, but the multipliers in question are described by the spherical transform which reduces to a commutative convolution subalgebra of functions on $G/K$.

\smallskip

Say that $T_m$ is an {\it $L^p-L^q$ Fourier multiplier} if it extends to a bounded map from $L^p( M)$ to $L^q( M)$.

\begin{theorem}[H. Zhang \cite{haonan.zhang}]\label{zhang-mult}
For $1< p\le 2 \le q<\infty$ and $r$ such that $1/r = 1/p-1/q$,
every $m\in L^{r,\infty}( \hat M )$ is an $L^p-L^q$ Fourier multiplier with
$$
\| T_m\|_{ L^p(M) \to L^q(M)} \le C_{p,q}\, \| m\|_{ L^{r,\infty}( \hat M ) }.
$$
\end{theorem}

\smallskip
\emph{\textbf{Example:}}
The case of $M=L^\infty(G)$ with a unimodular group $G$ was considered before by Akylzhanov and Ruzhansky \cite{ak-ruzh}.
They give also the following example: on the Heisenberg group $\mathbb H$ with the sub-Laplacian $\cal L = X^2 + Y^2$,
$$
\| (1-\cal L) ^{-s} \|_{ L^{r,\infty}( VN(\mathbb H) ) } \lesssim \sup_{u>0} u^{3/2r} (1+u)^{-s},
$$
which is finite as soon as $s\ge \dfrac 3{2r}$, so that in these assumptions $(1-\cal L) ^{-s}$ is bounded from $L^p(\mathbb H)$ to $L^q(\mathbb H)$.
The estimate of the weak $L^r$ norm is based on \eqref{tau_L}.

\smallskip
In the dual case of $M=VN(G)$ on a unimodular group $G$, the symbol of a multiplier is a bounded function $m:G\to\C$. The map
\begin{equation}\label{Tm}
T_m: L_f \mapsto L_{mf}
\end{equation}
has the above form: here $\cal F: L_f\mapsto f$ is the Fourier transform of $M$, which is the inverse to that of $\hat M = L^\infty(G)$, so that $L_{mf} = \cal F^{-1}( m \cal F(f))$.

\smallskip
\emph{\textbf{Free group of $n$ generators $G=\mathbb F_n$}}

Let $|g|$ be the length of the reduced word $g\in G$, and
$$
m_t(g) = e^{ -t|g|}.
$$
For positivity reasons, $T_{m_t}$ is bounded on $L^p(VN(G))$ for every $p\in[1,+\infty]$.
We can calculate the weak $L^r$ norm of $m_t$.

\medskip
We have $|m_t(g)| \ge\a$ %iff $-t|g| \ge \log\a$ 
iff $|g| \le - \frac1t \log \a = \frac1t | \log \a|$ since $m_t<1$.\\
The next formula is valid for $0<\a<1$, otherwise $\l_m(\a)=0$:
%$\# \{g: |m_t(g) \ge \a\} \le (2n)^{ - t^{-1} \log \a} $
$$
\l_m(\a) = |\{ g\in G: |m(g)|>\a \}| \le (2n)^{  t^{-1} |\log \a|}.
$$
$$
\| m_t\| _{ L^{r,\infty} (L^\infty(G)) } = \sup_{\a>0} \a \l_m^{1/r}(\a)
\le \sup_{0<\a<1} \a^{1- (tr)^{-1} \log (2n) }
$$
and this is finite iff $1- (tr)^{-1} \log (2n)\ge0$, that is, iff
$
t \ge \frac {\log (2n)} r
$.
Under this condition and $p,q$ as in Theorem \ref{zhang-mult}, $T_m$ is an $L^p-L^q$ multiplier.

Even more is true, on any discrete group (H. Zhang, \cite{haonan.zhang}): if $|\xi|\le m_t$ and
$1<p<\infty$, $\frac1{p*} = |\frac12-\frac1p|$, then
$\| T_\xi\| _{L^p( VN(G)) } \lesssim
\| \xi\|_{L^{p*,\infty}(G)} \le \| m_t\|_{L^{p*,\infty}(G)} $.

By different methods, another estimate was obtained by Junge, Palazuelos, Parcet, and Perrin \cite{JPPP}:
For $1< p\le q< \infty$, $ \| T_{m_t} \| _{ L^p(VN(G)) \to L^q( VN(G) )} \le 1$ as soon as $t> C \log \frac{ q-1}{p-1}$.

\smallskip
For \emph{\textbf{$L^p$ multipliers}}, a major result on $\R^d$ is
\begin{theorem}[H\"ormander--Mikhlin, 1956 -- 1960]
If $m$ is $\frac d2+1$ times differentiable and with some $C_m>0$, with usual notations for $\a\in \N^d$, %for $\a\in\N^n$, set $|\alpha|=\a_1+\dots+\a_n$.
\begin{equation}\label{HM}
\max_{0\le|\alpha|\le [\frac d2]+1} \sup_{\xi\in\R^d} |\xi|^{|\a|} |\partial^\a m(\xi)| \le C_m,
\end{equation}
then $m$ is an $L^p(\R^d)$ multiplier for every $1<p<\infty$.
\end{theorem}

Briefly, the estimate depends on a finite number of derivatives of $m$ which should decrease polynomially.
It is applicable in particular to $e^{-|\xi|}$ or to $(1+|\xi|^2)^{-r}$ with $r$ large enough. 
It turns out that similar assumptions have sense for multipliers on group algebras.

\smallskip
\emph{\textbf{Free group $\mathbb F_\infty$ of a countable set of generators}}

Let $\{g_k\}_{k\in\N}$ be the generators of $G$.
Let $m$ be a function on $\R^d$, satisfying \eqref{HM}.
For $g = g_{i_1}^{k_1} \dots g_{i_n}^{k_n}$, set
$$
\tilde m(g) = \begin{cases}
m( k_1, \dots, k_d) & \text{ if } n\le d;\\
m( k_1, \dots, k_n,0,\dots,0) & \text{ if } n<d. 
\end{cases}
$$
Then $T_{\tilde m}$ is bounded on $L^p( VN(\mathbb F_\infty))$ (Mei, Ricard, Xu \cite{xu-2022}).

\medskip
\emph{\textbf{Multipliers on the group algebra of $SL_n(\R)$}}

\let\g\gamma

For $g\in SL_n(\R)$, set
$
|g| := \max( \| g-e\|, \|g^{-1}-e\| )
$.
Let $X_1,\dots, X_{n^2-1}$ be the standard basis of
$$
\mathfrak{sl}_n(\R) = \{ A\in M_n(\R): Tr\, A=0\}.
$$
For $\g = ( j_1,\dots, j_{|\g|} )$ and a function $m$ denote
$$
d^\g_g m(g) = X_{j_1} \dots X_{j_|\g|} m.
$$

\begin{theorem}[Parcet, Ricard, De la Salle] \cite{PRS}
If for all $|\g| \le [\frac{n^2} 2]+1$ and $g\in G$
$$
 |g|^{|\g|} |d^\g_g m(g) | \le C_m,
$$
then $T_m$ is completely bounded on $L^p(VN(G))$ for $1<p<\infty$.
\end{theorem}

These examples close our brief survey.

\end{document}